\documentclass[12pt]{article} 
\usepackage[utf8]{inputenc}
\usepackage{amsmath,amssymb,amsthm,url,hyperref}

\theoremstyle{plain}
\newtheorem{theorem}{Theorem}

\theoremstyle{definition}

\theoremstyle{remark}
\newtheorem{remark}[theorem]{Remark}

\allowdisplaybreaks
\def\modd#1 #2{#1\ \mbox{\rm (mod}\ #2\mbox{\rm )}}

\begin{document}
\begin{center}
\vskip 1cm{\Large\bf 
Triangular Numbers With a Single Repeated Digit
}
\vskip 1cm
\large
Christian Hercher\\
Institut f\"{u}r Mathematik\\
Europa-Universit\"{a}t Flensburg\\
Auf dem Campus 1c\\
24943 Flensburg\\
Germany \\
\href{mailto:christian.hercher@uni-flensburg.de}{\tt christian.hercher@uni-flensburg.de} \\
Karl Fegert\\
Neu-Ulm\\
Germany\\
\href{mailto:karl.fegert@arcor.de}{\tt karl.fegert@arcor.de} \\
\end{center}

\begin{abstract}
The question of which triangular numbers have a decimal representation containing a single repeated digit seamed to be settled since at least the 1970s: Ballew and Weger provided a complete list and a proof that these are the only numbers of this kind. This assertion is referenced by other authors in the field. However, their proof is flawed. We present a new and elementary proof of the statement, which corrects the error.
\end{abstract}

\section{Introduction}

The study of triangular numbers is a long-standing and well-researched area within the field of number theory. For example, Gau\ss{} proved in his Disquisitiones Arithmeticae \cite{Gauss} that every positive integer can be expressed as the sum of three triangular numbers.

A further question that arises is that of the decimal representation of triangular numbers, in particular: Which triangular numbers can be expressed using a single digit? That is to say, are there triangular numbers of the form $d\cdot \frac{10^i-1}{9}$ with some digit $d\in\{1,\dots,9\}$ and postivie integer~$i$? One instance for this question is Problem~15648 \cite{Escott}, which was proposed by Youngman and subsequently answered by Escott in the same journal. However, Escott demonstrated only that triangular numbers with repeated digits exceeding $666$, comprising a maximum of $30$ digits, do not exist.

Subsequently, Ballew and Weger \cite{Ballew}  presented an erroneous proof for the conjecture that the only triangular numbers whose decimal representations consist of a single repeated digit are $1, 3, 6, 55, 66,$ and $666$. 

In a recent contribution,  Kafle, Luca, and Togb\'e \cite{Kafle} generalized the problem to include repeating blocks with two digits: The only triangular numbers whose decimal representation can be written using only one block of two digits are $10, 15, 21, 28, 36, 45, 55, 66, 78, 91, 5050,$ and $5151$. However, for the original problem, they just cite the work of Ballew and Weger.

Recently, on the internet platform MathOverflow \cite{MO} a discussion regarding this topic was held. In this discussion, the error in Ballew and Weger's proof was identified, and a proof using the methods of elliptic curves was presented. In this paper, we provide an elementary proof for this statement using only Pell equations.

\subsection*{Outline}
In accordance with the methodology proposed by Escott, Ballew, and Weger, we initially reformulate the problem. Subsequently, we provide a concise account of the erroneous assertion present within the proof of Ballew and Weger. In the following section, a proof is provided for each remaining digit, demonstrating that there are no further triangular numbers consisting of a single repeated digit.
 
\section{Reduction and reformulation of the problem}
For a positive integer $k$ let $T_k$ denote the $k$'th triangular number. In accordance with Escott's approach and the valid part of Ballew's and Weger's paper, we arrive at the following conclusion: If the decimal expression of $T_k$ consists solely of the single repeated digit $d$, we have
\begin{align*}
T_k=\frac{k(k+1)}{2}&=d\cdot \frac{10^i-1}{9}, i\geq 1, d\in\{1,\dots,9\}\\
\iff (2k+1)^2=4k^2+4k+1&=8d\cdot \frac{10^i-1}{9}+1.
\end{align*}

A solution $(k,d,i)$ in positive integers for the original problem exists, if and only if the quantity $D:=1+8d\cdot \frac{10^i-1}{9}$ is a perfect square.

It is $D \equiv \modd{1+8d} {10}$. Thus, immediately, we can rule out $d \in\{ 2, 4, 7, 9\}$, because in those cases $D$ would be congruent $7$ or $3$ modulo $10$, which are not quadratic residues. Therefore, the value of $d$ must be $1, 3, 5, 6,$ or $8$. Now we look at these cases one by one:

\begin{itemize}
\item[$d = 1$:] If $i=1$ it follows $D=8 \cdot \frac{10^i-1}{9}+1=9=3^2$, which yields $T_1 = 1$, and for $i>1$ we have $D = 88\ldots89$.

\item[$d = 3$:] If $i=1$ then $D=24 \cdot \frac{10^i-1}{9}+1=25=5^2$, which yields $T_2 = 3$. Otherwise it is $D = 26\ldots65$.

\item[$d = 5$:] If $i=1$ we have $D=40 \cdot \frac{10^i-1}{9}+1=41$. For $i=2$ we get $D=441=21^2$ ,which yields $T_{10} = 55$, and for $i>2$ it is $D = 44\ldots41$.

\item[$d = 6$:] If $i=1$, $2$, or $3$ then $D=48 \cdot \frac{10^i-1}{9}+1=49=7^2$, $D = 529 = 23^2$, or $D = 5329 = 73^2$, which yields $T_3 = 6, T_{11} = 66$, and $T_{36} = 666$. If $i>3$ it is $D = 53\ldots329$.

\item[$d = 8$:] If $i=1$ or $2$ then $D=64 \cdot \frac{10^i-1}{9}+1=65$ or $D=705$. Otherwise $D = 71\ldots105$. 
\end{itemize}

Since $\ldots05$ and $\ldots65$ cannot be the last two digits of a square number, it can be concluded that $d = 8$ can be excluded and $d = 3$ as well, except for $T_3$.

\subsection*{The false statement}
Ballew and Weger have incorrectly asserted that there is no number $z = \ldots88889$ whose square has more than four final digits $8$. This assertion is, however, erroneous. For example, $8072917^2 = \ldots88888889$. In fact, since there are solutions to the congruence equations $p^2 \equiv \modd{8 \cdot \frac{-1}{9} + 1} {2^3}$ and $p^2 \equiv \modd{8 \cdot \frac{-1}{9} + 1} {5^1}$, Hensel's Lemma ensures solutions modulo every power of $2$ and $5$, and thus, modulo every power of $10$. Conequently, there are square numbers whose last digits are $\ldots 8889$ with an arbitrary number of digits~$8$.

\section{A proof using Pell equations}
We therefore introduce a new line of reasoning which rules out the remaining cases, namely those with $d \in \{1,5,6\}$.

\subsection*{The case $d=1$}
\[p^2=88\ldots89=8\cdot \frac{10^i-1}{9}+1 \iff (3p)^2-(8\cdot 10^i)=1\]
Assume there is a solution $(p,i)$, $i \geq 2$ (the case $i = 1$ leads to $T_1$): Then we can make a case distinction by the parity of $i$.
\begin{itemize}
\item[Case A:] $i$ even, i.e., $i = 2r$, $r\in\mathbb{N}$. Subsequently, we have that 
\begin{align}
(3p)^2-2\cdot(2\cdot10^r)^2 &= 1, \notag 
\intertext{and we seek for solutions to the Pell equation}
x^2-2y^2 &= 1 \label{eq_1}
\end{align}
 with the requirements $x = 3p$ and $y = 2\cdot 10^r, r \geq 1$. Thus, $y$ must satisfy $y \equiv \modd{0} {5}$ and $y \not\equiv \modd{0} {7}$.

But the solutions to \eqref{eq_1} are given by $(x_0,y_0) = (3,2)$, $x_{n+1} = 3x_n + 4y_n$, $y_{n+1} = 2x_n + 3y_n$. A calculation modulo 5 and modulo 7 as in Table~\ref{Tab_1A} yields $y \equiv \modd{0} {5} \iff y \equiv \modd{0} {7}$, hence $y \neq 2\cdot 10^r$.

\begin{table}[hbt] 
\begin{center}
\begin{tabular}{c|cc|cc}
$n$ & $x_n\bmod{5}$ & $y_n\bmod{5}$ & $x_n\bmod{7}$ & $y_n\bmod{7}$\\ 
\hline
0 & 3 & 2 & 3 & 2\\
1 & 2 & 2 & 3 & 5\\
2 & 4 & 0 & 1 & 0\\
3 & 2 & 3 & 3 & 2\\
4 & 3 & 3 & 3 & 5\\
5 & 1 & 0 & 1 & 0\\
6 & 3 & 2 & 3 & 2\\
7 & \multicolumn{4}{c}{$\dots$}
\end{tabular}
\caption{Solutions of equation \eqref{eq_1} for $d=1$, Case A, mod $5$ and mod $7$.}
\label{Tab_1A}
\end{center}
\end{table}

\item[Case B:] $i$ odd, i.e., $i = 2r + 1$, $r\in\mathbb{N}$. Then we have $(3p)^2-20\cdot(2\cdot10^r)^2 = 1$, and we seek for solutions to 
\begin{align}
x^2 -20y^2 = 1 \label{eq_2}
\end{align}
 with $x = 3p$ and $y = 2  \cdot 10^r$, $r \geq 1$. Hence, it is $y \equiv \modd{0} {5}$ and $y \not\equiv \modd{0} {11}$.

But the solutions to \eqref{eq_2} are $(x_0,y_0) = (9,2)$, $x_{n+1} = 9x_n + 40y_n$, $y_{n+1} = 2x_n + 9y_n$, and a calculation modulo $5$ and $11$ similar to the one in Case~A yields $y \equiv \modd{0} {5} \iff y \equiv \modd{0} {11}$, so $y \neq 2 \cdot 10^r$ with $r \geq 1$.
\end{itemize}

\subsection*{The case $d=5$}
\[p^2=44\ldots41=40\cdot \frac{10^i-1}{9}+1.\]
Assume, there is a solution $(p,i)$.

\begin{itemize}
\item[Case A:] $i$ even, i.e., $i = 2r$, $r\in\mathbb{N}$. Then we have $(3p)^2-10\cdot(2\cdot10^r)^2 = -31$, and we seek solutions to the equation 
\begin{align}
x^2 -10y^2 &= -31 \label{eq_3}
\end{align}
 with $x = 3p$ und $y = 2  \cdot 10^r$, $r \geq 1$. Hence, $y$ has to satisfy the congruences $y \equiv \modd{0} {8}$ and $y \not\equiv \modd{0} {7}$.

It is evident, that  $(x_0,y_0) \in\{(3,2); (63,20)\}$ are solutions to \eqref{eq_3}. All other solutions of \eqref{eq_3} can be obtained from these by $x_{n+1} = 19x_n + 60y_n$, $y_{n+1} = 6x_n + 19y_n$. A calculation modulo $8$ and $7$ as in the case $d=1$ yields $y \equiv \modd{0} {8} \iff y \equiv \modd{0} {7}$, so $y \neq 2 \cdot 10^r$.

\item[Case B:] $i$ odd, i.e., the number of digits four in the decimal representation of $p^2$ is odd. Then according to a basic rule
\[44\ldots41\equiv 1 - 4 + (4 - 4) + \dots + (4 - 4) \equiv -3 \equiv \modd{8} {11},\] 
which is not a quadratic residue modulo~$11$.

(This is Problem 2 from the first round of the 2024 Bundeswettbewerb Mathematik in Germany \cite{BWM}. When we tried to solve the related open question for Case~A, i.e., whether 1 and 441 are the only squares with an even number of digits four, we came across Ballews and Wegers flawed proof and ended up with this paper.)
\end{itemize}

\subsection*{The case $d=6$}
Assume, there is a solution $(p,i)$ with $i\geq 4$. (The cases with $i\leq 3$ yield the solutions $T_3$, $T_{11}$, and $T_{36}$.) Then we have
\[p^2=53\ldots329=48\cdot \frac{10^i-1}{9}+1 \iff 9p^2-3 \cdot (4^2 \cdot 10^i)=-39.\]

\begin{itemize}
\item[Case A:] $i$ even, i.e., $i=2r$. Thus, $3p^2-(4\cdot 10^r)^2=-13$ and we seek solutions to the Pell equation 
\begin{align}
x^2-3y^2&=13 \label{eq_4}
\end{align}
 with the additional requirement of $x=4\cdot 10^r$, $y=3p$, and $r\geq 2$. Hence, it has to be $x\equiv \modd{0} {50}$. We get all solutions of \eqref{eq_4} starting from $(x_0,y_0)=(4,1)$ or $(x_0,y_0)=(5,2)$ through the recursion $x_{n+1}=2x_n+3y_n$ and $y_{n+1}=x_n+2y_n$. A calculation modulo $50$ and $241$ as in the case $d=1$ yields $x_i \equiv \modd{0} {50} \iff x_i \equiv \modd{94} {241}$. Thus, in every solution we have $94 \equiv \modd{4 \cdot 10^r} {241}$ and therefore, $10^r \equiv \modd{144} {241}$. However, this is not the case for any positive integer $r$.

\item[Case B:] $i$ odd, i.e., $i = 2r + 1$. This leads to $(3p)^2 - 30 \cdot (4 \cdot 10^r)^2 = -39$ and therefore, the Pell equation 
\begin{align}
x^2-30y^2&=-39  \label{eq_5}
\end{align}
 with $x=3p$, $y=4 \cdot 10^r$, and $r\geq 2$. All solutions of \eqref{eq_5} can be obtained from the starting solutions $(x_0,y_0)=(9,2)$ or $(x_0,y_0)=(21,4)$ through the recursion $x_{n+1}=11x_n+60y_n$, $y_{n+1}=2x_n+11y_n$. 

Now assume there is a solution with $r\geq 4$. Then it follows $y=4\cdot 10^r \equiv \modd{0} {64}$. A calculation modulo $64$ and $31$ as in the case $d=1$ yields $y_i\equiv \modd{0} {64} \iff y_i\equiv \modd{3} {31}$ and, therefore $10^r \equiv \modd{24} {31}$, which is not the case for any positive integer $r$.

In the cases with $r\leq 3$, it can be verified that in both families of solutions of \eqref{eq_5} we have $y_3>4000$. Consequently, the only remaining candidates for additional solutions are $(x_1,y_1)$ and $(x_2,y_2)$ in both families. However, a brief calculation reveals that they are no solutions to the original problem.
\end{itemize}

This demonstrates that, in each case, there are no additional solutions beyond those initially presented in the paper, thereby establishing the proof's conclusion.

\begin{remark}
The method can be applied in the generalized version with repeated blocks of digits as well. In consideration of blocks of length $2$, as presented in \cite{Kafle}, we have to solve the equation
\begin{alignat*}{2}
&\frac{k(k+1)}{2}&=&c\cdot \frac{100^i-1}{99}, \text{ $c\in\{10,\dots,99\}$}\\
\iff&(2k+1)^2&=&8c\cdot \frac{100^i-1}{99}+1\\
\iff&(33\cdot (2k+1))^2&=&88c \cdot 100^i + 1089 - 88c\\
\iff&x^2-22c \cdot y^2&=&1089-88c, \text{ with } x:=33\cdot (2k+1), y:=2\cdot 10^i
\end{alignat*}
in positive integers. Now we can compute the set of all integer solutions $(x,y)$ to the Pell equation given in the last line and find moduli that disprove that $y$ can be of the form $y=2\cdot 10^i$.
\end{remark}


\begin{thebibliography}{WWW}
\bibitem{Ballew} D. W. Ballew and R. C. Weger, Triangular numbers with repeated digits, \textit{Proc. S.D. ACAD. Sci.} \textbf{51} (1972), 52--55, \url{https://sdaos.org/wp-content/uploads/pdfs/Vol\%2051\%201972/72p52.pdf}.
\bibitem{Escott} E. B. Escott, and C. E. Youngman, Problem~15648, \textit{Educational Times} \textbf{8} (1905), 33--34, \url{https://archive.org/details/educationaltimes58educ/page/87/mode/1up}.
\bibitem{BWM} K. Fegert, Bundeswettbewerb Mathematik -- Aufgaben und Lösungen 1. Runde 2024, \url{https://www.mathe-wettbewerbe.de/fileadmin/Mathe-Wettbewerbe/Bundeswettbewerb_Mathematik/Dokumente/Aufgaben_und_Loesungen_BWM/2024/loes_24_1_e.pdf}, 2024.
\bibitem{Gauss} C. F. Gau\ss, \textit{Disquisitiones Arithmeticae}, In commiss. apud Gerh. Fleischer, 1801.
\bibitem{Kafle}B. Kafle, F. Luca, and A. Togb\'e, Triangular repblocks, \textit{Fibonacci 
Quart.} \textbf{56} (2018), 325--328, \url{https://www.fq.math.ca/Papers1/56-4/kafle1022018.pdf}.
\bibitem{MO} Internet discussion on mathoverflow, triangular repdigits, \url{https://mathoverflow.net/questions/467405/triangular-repdigits}, 2024.
\end{thebibliography}
\end{document}